\documentclass{amsart}
\usepackage{amsopn}
\usepackage{amssymb, amscd}

\newcommand{\nc}{\newcommand}

\nc{\vg}{\mathfrak{v} } \nc{\wg}{\mathfrak{w} }
\nc{\zg}{\mathfrak{z} } \nc{\ngo}{\mathfrak{n} }
\nc{\kg}{\mathfrak{k} } \nc{\mg}{\mathfrak{m} }
\nc{\bg}{\mathfrak{b} } \nc{\ggo}{\mathfrak{g} }
\nc{\ggob}{\overline{\mathfrak{g}} } \nc{\sog}{\mathfrak{so} }
\nc{\sug}{\mathfrak{su} } \nc{\spg}{\mathfrak{sp} }
\nc{\slg}{\mathfrak{sl} } \nc{\glg}{\mathfrak{gl} }
\nc{\cg}{\mathfrak{c} } \nc{\rg}{\mathfrak{r} }
\nc{\hg}{\mathfrak{h} } \nc{\tg}{\mathfrak{t} }
\nc{\ug}{\mathfrak{u} } \nc{\dg}{\mathfrak{d} }
\nc{\ag}{\mathfrak{a} } \nc{\pg}{\mathfrak{p} }
\nc{\sg}{\mathfrak{s} } \nc{\pca}{\mathcal{P}}
\nc{\nca}{\mathcal{N}} \nc{\lca}{\mathcal{L}} \nc{\oca}{\mathcal{O}}
\nc{\mca}{\mathcal{M}} \nc{\tca}{\mathcal{T}} \nc{\aca}{\mathcal{A}}
\nc{\cca}{\mathcal{C}} \nc{\sca}{\mathcal{S}}

\nc{\vp}{\varphi} \nc{\ddt}{{\small \frac{{\rm d}}{{\rm d}t}}}
\nc{\im}{\mathtt{i}}

\nc{\SO}{{\mathrm SO}} \nc{\Spe}{{\mathrm Sp}} \nc{\Sl}{{\mathrm
SL}} \nc{\SU}{{\mathrm SU}} \nc{\Or}{{\mathrm O}} \nc{\U}{{\mathrm
U}} \nc{\Gl}{{\mathrm GL}} \nc{\Se}{{\mathrm S}} \nc{\Cl}{{\mathrm
Cl}} \nc{\Spein}{{\mathrm Spin}} \nc{\Pin}{{\mathrm Pin}}

\nc{\RR}{{\Bbb R}} \nc{\HH}{{\Bbb H}} \nc{\CC}{{\Bbb C}}
\nc{\ZZ}{{\Bbb Z}} \nc{\FF}{{\Bbb F}} \nc{\NN}{{\Bbb N}}
\nc{\QQ}{{\Bbb Q}} \nc{\PP}{{\Bbb P}}

\nc{\vs}{\vspace{.2cm}} \nc{\vsp}{\vspace{1cm}}
\nc{\ip}{\langle\cdot,\cdot\rangle} \nc{\la}{\langle}
\nc{\ra}{\rangle} \nc{\unm}{\frac{1}{2}} \nc{\unc}{\frac{1}{4}}
\nc{\und}{\frac{1}{16}} \nc{\no}{\vs\noindent}
\nc{\lam}{\Lambda^2\ngo^*\otimes\ngo} \nc{\tangz}{{\rm T}^{\rm Zar}}
\nc{\nor}{{\sf n}} \nc{\eigen}{(k_1<...<k_r;d_1,...,d_r)}
\nc{\eigencero}{(0<k_2<...<k_r;d_1,...,d_r)} \nc{\mum}{/\!\!/}
\nc{\kir}{/\!\!/\!\!/}

\nc{\He}{\operatorname{Hess}} \nc{\ad}{\operatorname{ad}}
\nc{\Ad}{\operatorname{Ad}} \nc{\rank}{\operatorname{rank}}
\nc{\Irr}{\operatorname{Irr}} \nc{\End}{\operatorname{End}}
\nc{\Aut}{\operatorname{Aut}} \nc{\Inn}{\operatorname{Inn}}
\nc{\Der}{\operatorname{Der}} \nc{\Ker}{\operatorname{Ker}}
\nc{\Iso}{\operatorname{I}} \nc{\Diff}{\operatorname{D}}
\nc{\Lie}{\operatorname{L}} \nc{\tr}{\operatorname{tr}}
\nc{\dif}{\operatorname{d}} \nc{\sen}{\operatorname{sen}}
\nc{\modu}{\operatorname{mod}} \nc{\Ric}{\operatorname{Ric}}
\nc{\Ricac}{\operatorname{Ric^{ac}}}
\nc{\Ricg}{\operatorname{Ric^{\gamma}}}
\nc{\Ricc}{\operatorname{Ric^{c}}} \nc{\sym}{\operatorname{sym}}
\nc{\symac}{\operatorname{sym^{ac}}}
\nc{\symc}{\operatorname{sym^{c}}} \nc{\scalar}{\operatorname{sc}}
\nc{\grad}{\operatorname{grad}} \nc{\ricci}{\operatorname{ric}}
\nc{\ricciac}{\operatorname{ric^{ac}}}
\nc{\riccic}{\operatorname{ric^{c}}}
\nc{\riccig}{\operatorname{ric^{\gamma}}}
\nc{\Rin}{\operatorname{M}} \nc{\Le}{\operatorname{L}}
\nc{\tang}{\operatorname{T}} \nc{\level}{\operatorname{level}}
\nc{\rad}{\operatorname{r}} \nc{\abel}{\operatorname{ab}}

\newtheorem{theorem}{Theorem}[section]
\newtheorem{proposition}[theorem]{Proposition}
\newtheorem{corollary}[theorem]{Corollary}
\newtheorem{lemma}[theorem]{Lemma}
\newtheorem{definition}[theorem]{Definition}
\newtheorem{remark}[theorem]{Remark}

\newtheorem{example}[theorem]{Example}

\title[A canonical compatible metric]{A canonical compatible metric for geometric structures
on nilmanifolds}

\author{Jorge Lauret}

\address{Department of Mathematics, Yale University,
10 Hillhouse Box 208283 New Haven, CT 06520 USA.}

\address{FaMAF and CIEM, Universidad Nacional de C\'ordoba, Haya de la Torre s/n, 5000 C\'ordoba, Argentina}
\email{lauret@mate.uncor.edu}

\thanks{
Supported by CONICET and Guggenheim Foundation fellowships, and a
grant from FONCyT (Argentina).}

\begin{document}

\maketitle

\begin{abstract}
Let $(N,\gamma)$ be a nilpotent Lie group endowed with an invariant
geometric structure (cf. symplectic, complex, hypercomplex or any of
their `almost' versions).  We define a left invariant Riemannian
metric on $N$ compatible with $\gamma$ to be {\it minimal}, if it
minimizes the norm of the invariant part of the Ricci tensor among
all compatible metrics with the same scalar curvature.  We prove
that minimal metrics (if any) are unique up to isometry and scaling,
they develop soliton solutions for the `invariant Ricci' flow and
are characterized as the critical points of a natural variational
problem. The uniqueness allows us to distinguish two geometric
structures with Riemannian data, giving rise to a great deal of
invariants.

Our approach proposes to vary Lie brackets rather than inner
products; our tool is the moment map for the action of a reductive
Lie group on the algebraic variety of all Lie algebras, which we
show to coincide in this setting with the Ricci operator.  This
gives us the possibility to use strong results from geometric
invariant theory.

\end{abstract}

\section{Introduction}\label{intro}

Invariant structures on nilpotent Lie groups, as well as on their
compact versions, nilmanifolds (i.e. quotients by cocompact discrete
subgroups), play an important role in symplectic and complex
geometry.  Such structures are described easily in terms of data
arising from the Lie algebra, but can in turn be used to get very
interesting and exotic examples such as, to mention one, metrics
with exceptional holonomy.  The aim of this paper is the search for
the `best' Riemannian metric compatible with a given geometric
structure. Symplectic, complex and hypercomplex cases will be
treated in some detail, and although we will always have these
particular cases in mind, the main results will be proved in a more
general setting.

\subsection{Geometric structures and compatible metrics}
Let $N$ be a real $n$-di\-men\-sio\-nal nilpotent Lie group with Lie
algebra $\ngo$, whose Lie bracket will be denoted by $\mu
:\ngo\times\ngo\mapsto\ngo$.  By an invariant geometric structure on
$N$ we mean a tensor on $N$ defined by left translation of a tensor
$\gamma$ on $\ngo$ (or a set of tensors), usually non-degenerate in
some way, which satisfies a suitable integrability condition
\begin{equation}\label{closedgint}
{\rm IC}(\gamma,\mu)=0,
\end{equation}
involving only $\mu$ and $\gamma$.  The pair $(N,\gamma)$ will be
called a {\it class-$\gamma$ nilpotent Lie group}, and $N$ will be
assumed to be simply connected for simplicity.  A left invariant
Riemannian metric on $N$ is said to be {\it compatible} with
$(N,\gamma)$ if the corresponding inner product $\ip$ on $\ngo$
satisfies an orthogonality condition
\begin{equation}\label{ortconint}
{\rm OC}(\gamma,\ip)=0,
\end{equation}
in which only $\ip$ and $\gamma$ are involved.  We denote by
$\cca=\cca(N,\gamma)$ the set of all left invariant metrics on $N$
which are compatible with $(N,\gamma)$.  The pair $(\gamma,\ip)$
with $\ip\in\cca$ will often be referred to as a {\it class-$\gamma$
metric structure}.

A natural question arises:
\begin{itemize}
\item[] Given a class-$\gamma$ nilpotent Lie group $(N,\gamma)$, are there
canonical or distinguished left invariant Riemannian metrics on $N$
compatible with $\gamma$?,
\end{itemize}
where the meaning of 'canonical' and `distinguished' is of course
part of the problem.  The Ricci tensor has always been a very useful
tool to deal with this kind of questions, and since the answer
should depend on the metric and on the structure under
consideration, we consider the {\it invariant Ricci operator}
$\Ricg_{\ip}$ (and the invariant Ricci tensor
$\riccig_{\ip}=\la\Ricg_{\ip}\cdot,\cdot\ra$), that is, the
orthogonal projection of the Ricci operator $\Ric_{\ip}$ onto the
subspace of those symmetric maps of $\ngo$ leaving $\gamma$
invariant. D. Blair, S. Ianus and A. Ledger \cite{BlrIns, BlrLdg}
have proved in the compact case that metrics satisfying
\begin{equation}\label{unconditionint}
\riccig_{\ip}=0
\end{equation}
are very special in symplectic (so called metrics with hermitian
Ricci tensor) and contact geometry, as they are precisely the
critical points of two very natural curvature functionals on $\cca$:
the total scalar curvature functional $S$ and a functional $K$
measuring how far are the metrics of being K$\ddot{{\rm a}}$hler or
Sasakian, respectively (see also \cite{Blr}).

We will show that for a non-abelian nilpotent Lie group, condition
(\ref{unconditionint}) cannot hold for the classes of structures we
have in mind, and hence it is natural to try to get as close as
possible to this unattainable goal.  In this light, a metric
$\ip\in\cca(N,\gamma)$ is called {\it minimal} if it minimizes the
functional $||\riccig_{\ip}||^2=\tr(\Ricg_{\ip})^2$ on the set of
all compatible metrics with the same scalar curvature.  It turns out
that minimal metrics are the elements in $\cca$ closest to satisfy
the `Einstein-like' condition $\riccig_{\ip}=c\ip$, $c\in\RR$.  We
may also try to improve the metric via the evolution flow
$$
\ddt \ip_t=\pm\riccig_{\ip_t},
$$
whose fixed points are precisely metrics satisfying
(\ref{unconditionint}).  In the symplectic case, this flow is called
the anticomplexified Ricci flow and has been recently studied by H-V
Le and G. Wang \cite{LeWng}.  Of particular significance are then
those metrics for which the solution to the normalized flow (under
which the scalar curvature is constant in time) remains isometric to
the initial metric.  Such special metrics will be called {\it
invariant Ricci solitons}.  The main result in this paper can be now
stated.

\begin{theorem}\label{equiv2gint}
Let $(N,\gamma)$ be a nilpotent Lie group endowed with an invariant
geometric structure $\gamma$ (non-necessarily integrable). Then the
following conditions on a left invariant Riemannian metric $\ip$
which is compatible with $(N,\gamma)$ are equivalent:
\begin{itemize}
\item[(i)] $\ip$ is minimal.

\item[(ii)] $\ip$ is an invariant Ricci soliton.

\item[(iii)] $\Ricg_{\ip}=cI+D$ for some $c\in\RR$, $D\in\Der(\ngo)$.
\end{itemize}
Moreover, there is at most one compatible left invariant metric on
$(N,\gamma)$ up to isometry  (and scaling) satisfying any of the
above conditions.
\end{theorem}

A major obstacle to classify geometric structures is the lack of
invariants.  The uniqueness result in the above theorem gives rise
to a useful tool to distinguish two geometric structures; indeed, if
they are isomorphic then their respective minimal compatible metrics
(if any) have to be isometric.  One therefore can eventually
distinguish geometric structures with Riemannian data, which
suddenly provides us with a great deal of invariants.  This will be
used in Section \ref{exa} to find explicit continuous families of
pairwise non-isomorphic geometric structures in low dimensions,
mainly by using only one Riemannian invariant: the eigenvalues of
the Ricci operator.

A weakness of this approach is however the existence problem; the
theorem does not even suggest when such a distinguished metric does
exist.  How special are the symplectic or (almost-) complex
structures admitting a minimal metric?.  So far, we know how to deal
with this `existence question' only by giving several explicit
examples (see Remark \ref{part2}), for which the neat `algebraic'
characterization (iii) is very useful. It turns out that in low
dimensions the structures in general tend to admit a minimal
compatible metric, and the only obstruction we know at this moment
is when the subgroup of $\Gl(n)$ preserving $\gamma$ is in $\Sl(n)$
(e.g. for symplectic structures) and the nilpotent Lie algebra does
not have any nonzero symmetric derivation.  In this way, the
characteristically nilpotent Lie algebras (i.e. $\Der(\ngo)$ is
nilpotent) admitting a symplectic structure recently found by D.
Burde in \cite{Brd} can not be endowed with a minimal compatible
metric.  At the moment, these are the only counterexamples we have
to the existence question, and the lowest dimension among them is
$8$.

\subsection{Variety of compatible metrics and the moment map}\label{varint}
A class-$\gamma$ metric structure on a nilpotent Lie group is
determined by a triple $(\mu,\gamma,\ip)$ of tensors on $\ngo$.  The
proof of Theorem \ref{equiv2gint} is based on an approach which
proposes to vary the Lie bracket $\mu$ rather than the inner product
$\ip$.

Let us consider as a parameter space for the set of all real
nilpotent Lie algebras of a given dimension $n$, the set
$$
\nca=\{\mu\in V:\mu\;\mbox{satisfies Jacobi and is nilpotent}\},
$$
where $\ngo$ is a fixed $n$-dimensional real vector space and
$V=\lam$ is the vector space of all skew-symmetric bilinear maps
from $\ngo\times\ngo$ to $\ngo$. Since the Jacobi identity and the
nilpotency condition are both determined by zeroes of polynomials,
$\nca$ is a real algebraic variety.  We fix a tensor $\gamma$ on
$\ngo$ (or a set of tensors), and denote by $G_{\gamma}$ the
subgroup of $\Gl(n)$ preserving $\gamma$. The reductive Lie group
$G_{\gamma}$ acts naturally on $V$ leaving $\nca$ invariant and also
the algebraic subset $\nca_{\gamma}\subset\nca$ given by
$$
\nca_{\gamma}=\{\mu\in\nca:{\rm IC}(\gamma,\mu)=0\},
$$
that is, those nilpotent Lie brackets for which $\gamma$ is
integrable (see (\ref{closedgint})).

For each $\mu\in\nca$, let $N_{\mu}$ denote the simply connected
nilpotent Lie group with Lie algebra $(\ngo,\mu)$.  Fix an inner
product $\ip$ on $\ngo$ compatible with $\gamma$, that is, such that
(\ref{ortconint}) holds.  We identify each $\mu\in\nca_{\gamma}$
with a class-$\gamma$ metric structure on a nilpotent Lie group
\begin{equation}\label{idegint}
\mu\longleftrightarrow   (N_{\mu},\gamma,\ip),
\end{equation}
where all the structures are defined by left invariant translation.
The orbit $G_{\gamma}.\mu$ parameterizes then all the left invariant
metrics which are compatible with $(N_{\mu},\gamma)$ and hence we
may view $\nca_{\gamma}$ as the space of all class-$\gamma$ metric
structures on nilpotent Lie groups of dimension $n$. Two metrics
$\mu,\lambda\in\nca_{\gamma}$ are isometric if and only if they live
in the same $K_{\gamma}$-orbit, where
$K_{\gamma}=G_{\gamma}\cap\Or(\ngo,\ip)$ is the maximal compact
subgroup of $G_{\gamma}$.

We now go back to our search for the best compatible metric.  It is
natural to consider the functional $F:\nca_{\gamma}\mapsto\RR$ given
by $F(\mu)=\tr(\Ricg_{\mu})^2$, which in some sense measures how far
the metric $\mu$ is from satisfying (\ref{unconditionint}).  The
critical points of $F/||\mu||^4$ on the projective algebraic variety
$\PP\nca_{\gamma}\subset\PP V$ (which is equivalent to normalize by
the scalar curvature since $\scalar(\mu)=-\unc ||\mu||^2$), may
therefore be considered compatible metrics of particular
significance.

A crucial fact of this approach is that the moment map
$m_{\gamma}:V\mapsto\pg_{\gamma}$ for the action of $G_{\gamma}$ on
$V$ is proved to be
$$
m_{\gamma}(\mu)=8\Ricg_{\mu}, \qquad \forall\;\mu\in\nca_{\gamma},
$$
where $\pg_{\gamma}$ is the space of symmetric maps of $(\ngo,\ip)$
leaving $\gamma$ invariant (i.e.
$\ggo_{\gamma}=\kg_{\gamma}\oplus\pg_{\gamma}$ is a Cartan
decomposition).  This allows us to use strong and well-known results
on the moment map due to F. Kirwan \cite{Krw1} and L. Ness
\cite{Nss}, and proved by A. Marian \cite{Mrn} in the real case (see
Section \ref{git} for an overview).  Indeed, since $F$ becomes a
scalar multiple of the square norm of the moment map, we obtain the
following

\begin{theorem}\label{equiv1gint}\cite{Mrn}
Let $F:\PP \nca_{\gamma}\mapsto\RR$ be defined by
$F([\mu])=\tr(\Ricg_{\mu})^2/||\mu||^4$.  Then for $\mu\in
\nca_{\gamma}$ the following conditions are equivalent:
\begin{itemize}
\item[(i)] $[\mu]$ is a critical point of $F$.

\item[(ii)] $F|_{G_{\gamma}.[\mu]}$ attains its minimum value at $[\mu]$.

\item[(iii)] $\Ricg_{\mu}=cI+D$ for some $c\in\RR$, $D\in\Der(\mu)$.
\end{itemize}
Moreover, all the other critical points of $F$ in the orbit
$G_{\gamma}.[\mu]$ lie in $K_{\gamma}.[\mu]$.
\end{theorem}

The equivalence between (i) and (iii) in Theorem \ref{equiv2gint},
as well as the uniqueness result, follow then almost directly from
the above theorem.  We note that Theorem \ref{equiv1gint} also gives
a variational method to find minimal compatible metrics, by
characterizing them as the critical points of a natural curvature
functional (see Example \ref{m26} for an explicit application).

Most of the results obtained in this paper are still valid for
general Lie groups, although some considerations have to be
carefully taken into account (see Remark \ref{liegroups}).

\subsection{Examples}\label{exaint}
We first prove that a symplectic non-abelian nilpotent Lie group
$(N,\omega)$ can never admit a compatible left invariant metric with
hermitian Ricci tensor.  We also exhibit a curve of pairwise
non-isomorphic symplectic structures on the $6$-dimensional
nilpotent Lie group denoted by $(0,0,12,13,14+23,24+15)$ in
\cite{Slm}. Also, a curve of pairwise non-isomorphic non-abelian
complex structures on the Iwasawa manifold is given.  The initial
point of such a curve is the bi-invariant complex structure, and
after a finite period of time it becomes a curve of non-abelian
complex structures on the group denoted by $(0,0,0,0,12,14+23)$ in
\cite{Slm}.

By using results due to I. Dotti and A. Fino \cite{DttFin0,
DttFin2}, we prove that any hypercomplex $8$-dimensional nilpotent
Lie group admits a minimal compatible metric.  We actually show that
the moduli space of all hypercomplex $8$-dimensional nilpotent Lie
groups up to isomorphism is $9$-dimensional, and the moduli space of
the abelian ones has dimension $5$.  We finally give a surface of
pairwise non-isomorphic non-abelian hypercomplex structures on
$\ggo_3$, the $8$-dimensional Lie algebra obtained as the direct sum
of an abelian factor and the  $7$-dimensional quaternionic
Heisenberg Lie algebra.

\begin{remark}\label{part2}
{\rm More evidence of the existence of minimal compatible metrics is
showed in \cite{praga}, including all $4$-dimensional symplectic
structures and another curve in dimension $6$, two curves of abelian
complex structures on the Iwasawa manifold and several continuous
families depending on various parameters of abelian and non-abelian
hypercomplex structures in dimension $8$.  It is also showed in
\cite{praga} that if one considers no structure (i.e. $\gamma=0$),
then the `moment map' approach proposed in this paper can be also
applied to the study of Einstein solvmanifolds, obtaining many of
the uniqueness and structure results proved by J. Heber in
\cite{Hbr}.

By taking advantage again of the interplay with invariant theory, we
describe in \cite{classi} the moduli space of all isomorphism
classes of geometric structures on nilpotent Lie groups of a given
class and dimension admitting a minimal compatible metric, as the
disjoint union of semi-algebraic varieties which are homeomorphic to
categorical quotients of suitable linear actions of reductive Lie
groups. Such special geometric structures can therefore be
distinguished by using invariant polynomials.}
\end{remark}

\section{Geometric structures and compatible metrics}\label{geometric}

Let $N$ be a real $n$-dimensional nilpotent Lie group with Lie
algebra $\ngo$, whose Lie bracket is denoted by $\mu
:\ngo\times\ngo\mapsto\ngo$.  An invariant geometric structure on
$N$ is defined by left translation of a tensor $\gamma$ on $\ngo$
(or a set of tensors), usually non-degenerate in some way, which
satisfies a suitable integrability condition
\begin{equation}\label{closedg}
{\rm IC}(\gamma,\mu)=0,
\end{equation}
involving only $\mu$ and $\gamma$.  In this paper, we will focus on
the following classes of geometric structures: symplectic, complex
and hypercomplex, as well as on their respective `almost' versions,
that is, when condition (\ref{closedg}) is not required.  In this
way, ${\rm IC}(\gamma,\mu)$ can be for instance the differential of
a $2$-form or the Nijenhuis tensor associated to some
$(1,1)$-tensor. The contact case is somewhat different because the
condition is `open', but it becomes an equation of the form
(\ref{closedg}) when one considers fixed the underlying
almost-contact structure.  We shall deal with contact and complex
symplectic structures in a forthcoming paper.

The pair $(N,\gamma)$ will often be called a {\it class-$\gamma$
nilpotent Lie group}, and $N$ will be assumed to be simply connected
for simplicity.  The group $\Gl(n):=\Gl(n,\RR)=\Gl(\ngo)$ of
invertible maps of $\ngo$ acts on the vector space of tensors on
$\ngo$ of a given class, preserving the non-degeneracy, and if
$\gamma$ is integrable then $\vp.\gamma$ is so for any
$\vp\in\Aut(\ngo)$, the group of automorphisms of $\ngo$. In view of
this fact, two class-$\gamma$ nilpotent Lie groups $(N,\gamma)$ and
$(N',\gamma')$ are said to be {\it isomorphic} if there exists a Lie
algebra isomorphism $\vp:\ngo\mapsto\ngo'$ such that
$\gamma'=\vp.\gamma$. Also, given two geometric structures
$\gamma,\gamma'$ of the same class on $N$, we say that $\gamma$ {\it
degenerates to} $\gamma'$ if
$\gamma'\in\overline{\Aut(\ngo).\gamma}$, the closure of the orbit
$\Aut(\ngo).\gamma$ relative to the natural topology.

A left invariant Riemannian metric on $N$ is said to be {\it
compatible} with $(N,\gamma)$ if the corresponding inner product
$\ip$ on $\ngo$ satisfies an orthogonality condition
\begin{equation}\label{ortcon}
{\rm OC}(\gamma,\ip)=0,
\end{equation}
in which only $\ip$ and $\gamma$ are involved.  We denote by
$\cca=\cca(N,\gamma)$ the set of all left invariant metrics on $N$
which are compatible with $(N,\gamma)$.  The pair $(\gamma,\ip)$
with $\ip\in\cca$ will often be referred to as a {\it class-$\gamma$
metric structure}.  It is clear from (\ref{ortcon}) that for an
invariant geometric structure there always exist a compatible
metric, since the condition is independent from $\mu$. Moreover, the
space $\cca$ is usually huge; recall for instance that the group
$$
G_{\gamma}=\{ \vp\in\Gl(n):\vp.\gamma=\gamma\}
$$
acts on $\cca$, and it is easy to see that actually for any
$\ip\in\cca$ we have that
\begin{equation}\label{compg}
\cca=G_{\gamma}.\ip=\{\la\vp^{-1}\cdot,\vp^{-1}\cdot\ra:\vp\in
G_{\gamma}\}.
\end{equation}
A natural question takes place:
\begin{itemize}
\item[] Given a class-$\gamma$ nilpotent Lie group $(N,\gamma)$, are there
canonical or distinguished left invariant Riemannian metrics on $N$
compatible with $\gamma$?
\end{itemize}
This problem may be (and it is) stated for differentiable manifolds
in general, and does not only present some interest in Riemannian
geometry; indeed, the existence of a certain nice compatible metric
could eventually help to distinguish two geometric structures as
well as to find privileged geometric structures on a given manifold.

The aim of this section is to propose two properties which make a
compatible metric very distinguished, one is obtained by minimizing
a curvature functional and the other as a soliton solution for a
natural evolution flow. The Ricci tensor will be used in both
approaches.  In Appendix \ref{app}, we have reviewed some well known
properties of left invariant metrics on nilpotent Lie groups, which
will be used constantly from now on.

Fix a class-$\gamma$ nilpotent Lie group $(N,\gamma)$.  Let
$\ggo_{\gamma}$ be the Lie algebra of $G_{\gamma}$,
$$
\ggo_{\gamma}=\{ A\in\glg(n):A.\gamma=0\}.
$$

\begin{definition}\label{invricci}
{\rm For each compatible metric, we consider the orthogonal
projection $\Ricg_{\ip}$ of the Ricci operator $\Ric_{\ip}$ on
$\ggo_{\gamma}$,  called the {\it invariant Ricci operator}, and the
corresponding {\it invariant Ricci tensor} given by
$\riccig=\la\Ricg\cdot,\cdot\ra$. }
\end{definition}

The role of the Ricci tensor has always been crucial in defining
privileged (compatible) metrics; we have for example Einstein
metrics, extremal K$\ddot{{\rm a}}$hler metrics in complex geometry,
and more recently metrics with hermitian Ricci tensor (i.e. when the
Ricci operator commutes with $J$) and $\phi$-invariant Ricci tensor
in symplectic and contact geometry, respectively. These two last
notions are equivalent to $\riccig=0$ and have been characterized in
the compact case by D. Blair, S. Ianus and A. Ledger \cite{BlrIns,
BlrLdg} as the critical points of two very natural curvature
functionals on $\cca$: the total scalar curvature functional $S$ and
a functional $K$ for which the global minima are precisely
K$\ddot{{\rm a}}$hler or Sasakian metrics, respectively (see also
\cite{Blr}).

In this light, condition
\begin{equation}\label{uncondition}
\riccig_{\ip}=0,
\end{equation}
involves both the geometric structure and the metric, and seems to
be very natural to require to a compatible metric. Nevertheless, if
$\RR I\subset\ggo_{\gamma}$, then $\tr{\Ricg_{\ip}}=\scalar(\ip)$,
and so it is forbidden for instance for non-abelian nilpotent Lie
groups (where always $\scalar(\ip)<0$) in the complex and
hypercomplex cases. We shall prove that this condition is forbidden
in the symplectic case as well.  We therefore have to consider
(\ref{uncondition}) as an unreachable goal and try to get as close
as possible, for instance, by minimizing
$||\riccig_{\ip}||^2=\tr(\Ricg_{\ip})^2$. In order to avoid
homothetical changes, we must normalize the metrics some way.  In
the noncompact homogeneous case, the scalar curvature always appears
as a very natural choice.  We then propose the following

\begin{definition}\label{minimalg}
{\rm A left invariant metric $\ip$ compatible with a class-$\gamma$
nilpotent Lie group $(N,\gamma)$ is called {\it minimal} if
$$
\tr(\Ricg_{\ip})^2=\min\{  \tr(\Ricg_{\ip'})^2 :
\ip'\in\cca(N,\gamma), \quad \scalar(\ip')=\scalar(\ip)\}.
$$ }
\end{definition}

Recall that the existence and uniqueness (up to isometry and
scaling) of minimal metrics is far to be clear from the definition.
The uniqueness shall be proved in Section \ref{var}, but the
`existence question' is still nebulous.  Minimal metrics are the
compatible metrics closest to satisfy the `Einstein-like' condition
$\riccig_{\ip}=c\ip$, for some $c\in\RR$.   Indeed,
$$
||\Ricg_{\ip}-\frac{\tr{\Ricg_{\ip}}}{n}I||^2=\tr(\Ricg_{\ip})^2-\frac{(\tr{\Ricg_{\ip}})^2}{n}
$$
and $\tr{\Ricg_{\ip}}$ equals either $0$ or $\scalar(\ip)$,
depending on $\ggo_{\gamma}$ contains or not $\RR I$.

We now consider an evolution approach.  Motivated by the famous
Ricci flow introduced by R. Hamilton \cite{Hml1}, we consider the
{\it invariant Ricci flow} for our left invariant metrics on $N$,
given by the following evolution equation
\begin{equation}\label{grfn}
\ddt \ip_t=\pm\riccig_{\ip_t},
\end{equation}
which coincides for example with the anticomplexified  Ricci flow
studied in \cite{LeWng} in the symplectic case. The choice of the
best sign might depend on the class of structure.  This is just an
ordinary differential equation and hence the existence for almost
all $t$ and uniqueness of the solution is guaranteed. It follows
from (\ref{compg}) that
\begin{equation}\label{tangcompg}
\tang_{\ip}\cca=\{\alpha\in\sym(\ngo):A_{\alpha}.\gamma=0\},
\end{equation}
and therefore, if $\ip_0\in\cca$ then the solution $\ip_t\in\cca$
for all $t$ since $\riccig_{\ip_t}\in\tang_{\ip_t}\cca$ (see
Appendix \ref{app}).

\begin{remark}\label{uniq}
{\rm If we had however the uniqueness of the solution for the flow
(\ref{grfn}) in the non-compact general case, then we would not need
to restrict ourselves to left invariant metrics. Indeed, if $f$ is
an isometry of the initial metric $\ip_0$ which also leaves $\gamma$
invariant, then since $f^*\ip_t$ is also a solution and
$f^*\ip_0=\ip_0$ we would get by uniqueness of the solution that $f$
is an isometry of all the metrics $\ip_t$ as well. Left invariance
of the starting metric would be therefore preserved along the flow.
}
\end{remark}

When $M$ is compact, a normalized Ricci flow is often considered,
under which the volume of the solution metric is constant in time.
Actually, the normalized equation differs from the Ricci flow only
by a change of scale in space and a change of parametrization in
time (see \cite{Hml2, CaoChw}).  In our case, where the manifold is
non-compact but the metrics are homogeneous, it seems natural to do
the same thing but normalizing by the scalar curvature, which is
just a single number associated to the metric.  We recall that a
left invariant metric $\ip$ on a nilpotent Lie group $N$ has always
$\scalar(\ip)<0$, unless $N$ is abelian (see (\ref{ricci})).

\begin{proposition}\label{nacrf}
The solution to the normalized invariant Ricci flow
\begin{equation}\label{nfg}
\ddt\ip_t=\pm\riccig_{\ip_t}\mp\frac{\tr(\Ricg_{\ip_t})^2}{\scalar(\ip_t)}\ip_t
\end{equation}
satisfies $\scalar(\ip_t)=\scalar(\ip_0)$ for all $t$.  Moreover,
this flow differs from the invariant Ricci flow {\rm (\ref{grfn})}
only by a change of scale in space and a change of parametrization
in time.
\end{proposition}

\begin{proof}
It follows from the formula for the gradient of the scalar curvature
functional $\scalar:\pca\mapsto\RR$ given in (\ref{gradsc}) that if
$\ip_t$ is a solution of (\ref{nfg}), then the function
$f(t)=\scalar(\ip_t)$ satisfies
$$
\begin{array}{rl}
f'(t)&=g_{\ip_t}(\ddt\ip_t,-\ricci_{\ip_t})\\ \\
&=\mp
g_{\ip_t}(\riccig_{\ip_t},\ricci_{\ip_t})\pm\frac{\tr(\Ricg_{\ip_t})^2}{\scalar(\ip_t)}
g_{\ip_t}(\ip_t,\ricci_{\ip_t})\\ \\
&=\mp\tr(\Ricg_{\ip_t}\Ric_{\ip_t})\pm\frac{\tr(\Ricg_{\ip_t})^2}{\scalar(\ip_t)}\tr(\Ric_{\ip_t})\\ \\
&=\mp\tr(\Ricg_{\ip_t})^2(1-\frac{f(t)}{\scalar(\ip_t)})=0, \qquad
\forall t,
\end{array}
$$
and thus $f(t)\equiv f(0)=\scalar(\ip_0)$.  The last assertion
follows as in \cite{Hml2} in a completely analogous way.
\end{proof}

The fixed points of this normalized flow (\ref{nfg}) are those
metrics satisfying $\Ricg_{\ip}\in\RR I$, and so in particular, if
$\ggo_{\gamma}\subset\slg(n)$, then this is equivalent to
$\Ricg_{\ip}=0$.  Indeed, if
$\Ricg_{\ip}=\mp\frac{\tr(\Ricg_{\ip_t})^2}{\scalar(\ip_0)}I$ then
$\Ricg_{\ip}=0$ since $\tr\Ricg_{\ip}=0$.  We should also note that
for the flow (\ref{nfg}), $\ddt\ip_t\in\tang_{\ip_t}\cca+\RR I$ for
all $t$, which implies that the solution $\ip_t$ stays in the set of
all scalar multiples of compatible metrics.  Recall that if $\RR
I\subset\ggo_{\gamma}$ then the solution stays anyway in $\cca$.

In these evolution approaches always appear naturally the soliton
metrics, which are not fixed points of the flow but are close to,
and they play an important role in the study of singularities (see
the surveys \cite{Hml2, CaoChw} for further information).  The idea
is that if one is trying to improve a metric via an evolution
equation, then those metrics for which the solution remains
isometric to the initial point may be certainly considered as very
distinguished.

\begin{definition}\label{solitong}
{\rm A metric $\ip$ compatible with $(N,\gamma)$ is called an {\it
invariant Ricci soliton} if the solution $\ip_t$ to the normalized
invariant Ricci flow (\ref{nfg}) with initial metric $\ip_0=\ip$ is
given by $\vp_t^*\ip$, the pullback of $\ip$ by a one parameter
group of diffeomorphisms $\{\vp_t\}$ of $N$. }
\end{definition}

We now give a neat characterization of invariant Ricci soliton
metrics, which will be very useful in Section \ref{var} to prove the
equivalence with the property of being minimal (see Definition
\ref{minimalg}), and to find explicit examples in the subsequent
sections.

\begin{proposition}\label{soliton1g}
Let $(N,\gamma)$ be a class-$\gamma$ nilpotent Lie group.  A
compatible metric $\ip$ is an invariant Ricci soliton if and only if
$\Ricg_{\ip}=cI+D$ for some $c\in\RR$ and $D\in\Der(\ngo)$.  In such
a case, $c=\frac{\tr(\Ricg_{\ip})^2}{\scalar(\ip)}$.
\end{proposition}

\begin{proof}
We first note that the assertion on the value of the number $c$
follows from (\ref{ricort}); in fact,
$$
\tr(\Ricg_{\ip})^2=\tr(\Ric_{\ip}\Ricg_{\ip})=c\tr{\Ric_{\ip}}+\tr(\Ric_{\ip}D)=c\scalar(\ip).
$$
Assume that there exists a one-parameter group of diffeomorphisms
$\vp_t$ on $N$ such that $\ip_t=\vp^*_t\ip$ is a solution to the
flow (\ref{nfg}).  By the uniqueness of the solution we have that
$\vp_t^*g$ is also $N$-invariant for all $t$ (see Remark
\ref{uniq}).  Thus $\vp_t$ normalizes $N$ and so it follows from
\cite[Thm 2, 4)]{Wls} that $\vp_t\in\Aut(N).N$.  This implies that
there exists a one-parameter group $\psi_t$ of automorphisms of $N$
such that $\vp_t^*\ip=\psi_t^*\ip$ for all $t$.  Now, if
$\psi_t=e^{-\frac{t}{2}D}$ with $D\in\Der(\ngo)$ then
$\ddt|_0\psi_t^*\ip=\la D\cdot,\cdot\ra$, and using that
$\psi_t^*\ip$ is a solution of (\ref{nfg}) in $t=0$ we obtain that
$\riccig_{\ip}=c\ip+\la D\cdot,\cdot\ra$ for some $c\in\RR$, or
equivalently, $\Ricg=cI+D$.

Conversely, if $\Ricg=cI+D$ then we will show that the curve
$\ip_t=e^{-\frac{t}{2}D}.\ip$ is a solution of the flow (\ref{nfg}).
For any $t$, it follows from
$\frac{t}{2}D=\frac{t}{2}\Ricg_{\ip}-\frac{t}{2}cI\in\ggo_{\gamma}+\RR
I$ that
$$
\gamma=b(t)e^{-\frac{t}{2}D}.\gamma.
$$
for some $b(t)\in\RR$.  This implies that
$$
\Ricg_{\ip_t}=e^{-\frac{t}{2}D}\Ricg_{\ip}e^{\frac{t}{2}D}=e^{-\frac{t}{2}D}(cI+D)e^{\frac{t}{2}D}=cI+D
$$
for all $t$.  Therefore
$$
\begin{array}{rl}
\ddt|_0\ip_t&=\la D\cdot,\cdot\ra_t=\la(\Ricg_{\ip_t}-cI)\cdot,\cdot\ra_t \\ \\
&=\riccig_{\ip_t}-c\ip_t=\riccig_{\ip_t}+\frac{\tr(\Ricg_{\ip_t})^2}{\scalar(\ip_0)}\ip_t,
\end{array}
$$
as was to be shown.
\end{proof}

Recall that the condition in the above proposition can be replaced
by
$$
\Ricg_{\ip}-\frac{\tr(\Ricg_{\ip})^2}{\scalar(\ip)}I\in\Der(\ngo),
$$
which gives a computable method to check whether a metric is an
invariant Ricci soliton or not.

\section{Real geometric invariant theory and the moment map}\label{git}

In this section, we overview some results from (geometric) invariant
theory over the real numbers.  We refer to \cite{RchSld} for a
detailed exposition.  These will be our tools to study metrics
compatible with geometric structures on nilpotent Lie groups.

Let $G$ be a real reductive Lie group acting on a real vector space
$V$ (see \cite{RchSld} for a precise definition of the situation).
Let $\ggo$ denote the Lie algebra of $G$ with Cartan decomposition
$\ggo=\kg\oplus\pg$, where $\kg$ is the Lie algebra of a maximal
compact subgroup $K$ of $G$.  Endow $V$ with a fixed from now on
inner product $\ip$ such that $\kg$ and $\pg$ act by skew-symmetric
and symmetric transformations, respectively. Let $\mca=\mca(G,V)$
denote the set of {\it minimal vectors}, that is
$$
\mca=\{ v\in V: ||v||\leq ||g.v||\quad\forall g\in G\}.
$$
For each $v\in V$ define
$$
\rho_v:G\mapsto\RR, \qquad \rho_v(g)=||g.v||^2=\la g.v,g.v\ra.
$$
In \cite{RchSld}, R. Richardson and P. Slodowy showed that the nice
interplay between closed orbits and minimal vectors found by G.
Kempf and L. Ness for actions of complex reductive algebraic groups,
is still valid in the real situation.

\begin{theorem}\cite{RchSld}\label{RS}
Let $V$ be a real representation of a real reductive Lie group $G$,
and let $v\in V$.
\begin{itemize}
\item[(i)] $v\in\mca$ if and only if $\rho_v$ has a critical point at $e\in G$.

\item[(ii)] If $v\in\mca$ then $G.v\cap\mca=K.v$.

\item[(iii)] The orbit $G.v$ is closed if and only if $G.v$ meets $\mca$.
\end{itemize}
\end{theorem}

Let $(\dif\rho_v)_e:\ggo\mapsto\RR$ denote the differential of
$\rho_v$ at the identity $e$ of $G$. It follows from the
$K$-invariance of $\ip$ that $(\dif\rho_v)_e$ vanishes on $\kg$, and
so we can assume that $(\dif\rho_v)_e\in\pg^*$, the vector space of
real-valued functionals on $\pg$.  We therefore may define a
function called the {\it moment map} for the action of $G$ on $V$ by
\begin{equation}\label{moment}
m:V\mapsto\pg, \qquad  \la m(v),A\ra_{\pg}=(\dif\rho_v)_e(A),
\end{equation}
where $\ip_{\pg}$ is an $\Ad(K)$-invariant inner product on $\pg$.
Since $m(tv)=t^2m(v)$ for all $t\in\RR$, we also may consider the
moment map
\begin{equation}\label{moment2}
m:\PP V\mapsto\pg, \qquad  m(x)=\frac{m(v)}{||v||^2}, \quad 0\ne
v\in V,\;x=[v],
\end{equation}
where $\PP V$ is the projective space of lines in $V$.  If
$\pi:V\setminus\{ 0\}\mapsto \PP V$ denotes the usual projection
map, then $\pi(v)=x$.  In the complex case, under the natural
identifications $\pg=\pg^*=(\im\kg)^*=\kg^*$, the function $m$ is
precisely the moment map from symplectic geometry, corresponding to
the Hamiltonian action of $K$ on the symplectic manifold $\PP V$
(see for instance the survey \cite{Krw2} or \cite[Chapter 8]{Mmf}
for further information).

Consider the functional square norm of the moment map
\begin{equation}\label{norm}
F:V\mapsto\RR, \qquad  F(v)=||m(v)||^2=\la m(v),m(v)\ra_{\pg},
\end{equation}
which is easily seen to be a $4$-degree homogeneous polynomial.
Recall that $\mca$ coincides with the set of zeros of $F$. It then
follows from Theorem \ref{RS}, parts (i) and (iii), that an orbit
$G.v$ is closed if and only if $F(w)=0$ for some $w\in G.v$, and in
that case, the set of zeros of $F|_{G.v}$ coincides with $K.v$.  A
natural question arises: what is the role played by the remaining
critical points of $F:\PP V\mapsto\RR$ (i.e. those for which
$F(x)>0$) in the study of the $G$-orbit space of the action of $G$
on $V$, as well as on other real $G$-varieties contained in $V$?.
This was studied independently by F. Kirwan \cite{Krw1} and L. Ness
\cite{Nss}, and it is shown in the complex case that the non-minimal
critical points share some of the nice properties of minimal vectors
stated in Theorem \ref{RS}.  In the real case, which is actually the
one we need to apply in this paper, the analogous of some of these
results have been proved by A. Marian \cite{Mrn}.

 \begin{theorem}\cite{Mrn}\label{marian}
Let $V$ be a real representation of a real reductive Lie group $G$,
$m:\PP V\mapsto\pg$ the moment map and $F=||m||^2:\PP V\mapsto\RR$.
\begin{itemize}
\item[(i)] If $x\in\PP V$ is a critical point of $F$ then the functional $F|_{G.x}$ attains its minimum value at
$x$.

\item[(ii)] If nonempty, the critical set of $F|_{G.x}$ consists of a unique $K$-orbit.
\end{itemize}
\end{theorem}

We endow $\PP V$ with the Fubini-Study metric defined by $\ip$ and
denote by $x\mapsto A_x$ the vector field on $\PP V$ defined by
$A\in\ggo$ via the action of $G$ on $\PP V$, that is,
$A_x=\ddt|_0\exp(tA).x$.

\begin{lemma}\cite{Mrn}\label{marian2}
The gradient of the functional $F=\parallel m\parallel^2:\PP
V\mapsto\RR$ is given by
$$
\grad(F)_x=4m(x)_x, \qquad x\in\PP V,
$$
and therefore $x$ is a critical point of $F$ if and only if
$m(x)_x=0$, if and only if $\exp{tm(x)}$ fixes $x$.
\end{lemma}

We now develop some examples which are far to be the natural ones,
but they are those ones will be considered in this paper to study
left invariant structures on nilpotent Lie groups.

\begin{example}\label{gln}
{\rm Let $\ngo$ be an $n$-dimensional real vector space and $V=\lam$
the vector space of all skew-symmetric bilinear maps from
$\ngo\times\ngo$ to $\ngo$. There is a natural action of
$\Gl(n):=\Gl(n,\RR)$ on $V$ given by
\begin{equation}\label{action}
g.\mu(X,Y)=g\mu(g^{-1}X,g^{-1}Y), \qquad X,Y\in\ngo, \quad
g\in\Gl(n),\quad \mu\in V.
\end{equation}
Any inner product $\ip$ on $\ngo$ defines an $\Or(n)$-invariant
inner product on $V$, denoted also by $\ip$, as follows:
\begin{equation}\label{innV}
\la\mu,\lambda\ra=\sum_{ijk}\la\mu(X_i,X_j),X_k\ra\la\lambda(X_i,X_j),X_k\ra,
\end{equation}
where $\{ X_1,...,X_n\}$ is any orthonormal basis of $\ngo$.  A
Cartan decomposition for the Lie algebra of $\Gl(n)$ is given by
$\glg(n)=\sog(n)\oplus\sym(n)$, that is, in skew-symmetric and
symmetric transformations respectively, and we consider the
following $\Ad(\Or(n))$-invariant inner product on $\pg:=\sym(n)$,
\begin{equation}\label{inng}
\la A,B\ra_{\pg}=\tr{AB}, \qquad A,B\in\pg.
\end{equation}
The action of $\glg(n)$ on $V$ obtained by differentiation of
(\ref{action}) is given by
\begin{equation}\label{actiong}
A.\mu=-\delta_{\mu}(A):=
A\mu(\cdot,\cdot)-\mu(A\cdot,\cdot)-\mu(\cdot,A\cdot), \qquad
A\in\glg(n),\quad\mu\in V.
\end{equation}
If $\mu\in V$ satisfies the Jacobi condition, then
$\delta_{\mu}:\glg(n)\mapsto V$ coincides with the cohomology
coboundary operator of the Lie algebra $(\ngo,\mu)$ from level $1$
to $2$, relative to cohomology with values in the adjoint
representation.  Recall that $\Ker{\delta_{\mu}}=\Der(\mu)$, the Lie
algebra of derivations of the algebra $\mu$.  Let $A^t$ denote the
transpose relative to $\ip$ of a linear transformation
$A:\ngo\mapsto\ngo$ and consider the adjoint map
$\ad_{\mu}{X}:\ngo\mapsto\ngo$ (or left multiplication) defined by
$\ad_{\mu}{X}(Y)=\mu(X,Y)$.

\begin{proposition}\label{ourmoment}
The moment map $m:V\mapsto\pg$ for the action {\rm (\ref{action})}
of $\Gl(n)$ on $V=\lam$ is given by
\begin{equation}\label{defm}
m(\mu)=-4\displaystyle{\sum_{i}}(\ad_{\mu}{X_i})^t\ad_{\mu}{X_i}
+2\displaystyle{\sum_{i}}\ad_{\mu}{X_i}(\ad_{\mu}{X_i})^t,
\end{equation}
where $\{ X_1,...,X_n\}$ is any orthonormal basis of $\ngo$, and it
is a simple calculation to see that
\begin{equation}\label{defm2}
\begin{array}{rl}
\la m(\mu)X,Y\ra=&-4\displaystyle{\sum_{ij}}\la\mu(X,X_i),X_j\ra\la\mu(Y,X_i),X_j\ra \\
&+2\displaystyle{\sum_{ij}}\la\mu(X_i,X_j),X\ra\la\mu(X_i,X_j),
Y\ra, \qquad \forall\; X,Y\in\ngo.
\end{array}
\end{equation}
\end{proposition}

\begin{proof}
For any $A\in\pg$ we have that
$$
\begin{array}{rl}
(\dif\rho_{\mu})_I(A)=&\ddt|_0\la e^{tA}.\mu,e^{tA}.\mu\ra =-2\la\mu,\delta_{\mu}(A)\ra \\ \\

=&-2\displaystyle{\sum_{pij}}
\la\mu(X_p,X_i),X_j\ra\la\delta_{\mu}(A)(X_p,X_i),X_j\ra \\

=&-2\Big(\displaystyle{\sum_{pij}}\la\mu(X_p,X_i),X_j\ra\la\mu(AX_p,X_i),X_j\ra \\

&+ \la\mu(X_p,X_i),X_j\ra\la\mu(X_p,AX_i),X_j\ra  \\ \\

&-\la\mu(X_p,X_i),X_j\ra\la A\mu(X_p,X_i),X_j\ra\Big) \\ \\

=&
-2\Big(\displaystyle{\sum_{pijr}}\la\mu(X_p,X_i),X_j\ra\la\mu(X_r,X_i),X_j\ra\la AX_p,X_r\ra \\

& +\la\mu(X_p,X_i),X_j\ra\la\mu(X_p,X_r),X_j\ra\la AX_i,X_r\ra  \\ \\

& -\la\mu(X_p,X_i),X_j\ra\la\mu(X_p,X_i),X_r\ra\la AX_j,X_r\ra\Big).
\end{array}
$$
By interchanging the indexes $p$ and $i$ in the second line, and $p$
and $j$ in the third one, we get
$$
\begin{array}{rl}
(\dif\rho_{\mu})_I(A)=&
-4\displaystyle{\sum_{prij}}\la\mu(X_p,X_i),X_j\ra
\la\mu(X_r,X_i),X_j\ra
\la AX_p,X_r\ra\\

&+2\displaystyle{\sum_{prij}}\la\mu(X_i,X_j),X_p\ra\la\mu(X_i,X_j),X_r\ra\la
AX_p,X_r\ra.
\end{array}
$$
If we call $M$ the right hand side of (\ref{defm}), then we obtain
from (\ref{defm2}) that
$$
(\dif\rho_{\mu})_I(A)=\displaystyle{\sum_{pr}} \la MX_p,X_r\ra\la
AX_p,X_r\ra=\tr{MA}=\la M,A\ra_{\pg},
$$
which implies that $m(\mu)=M$ (see (\ref{moment})).
\end{proof}   }
\end{example}

\begin{example}\label{momg}
{\rm We keep the notation as in Example \ref{gln}.  Let $\gamma$ be
a tensor on $\ngo$ and let $G_{\gamma}\subset\Gl(n)$ denote the
subgroup leaving $\gamma$ invariant, with Lie algebra
$\ggo_{\gamma}$.  The group $G_{\gamma}$ is reductive and
$K_{\gamma}=G_{\gamma}\cap\Or(n)$ is a maximal compact subgroup of
$G_{\gamma}$, whose Lie algebra will be denoted by $\kg_{\gamma}$.
A Cartan decomposition is given by
$$
\ggo_{\gamma}=\kg_{\gamma}\oplus\pg_{\gamma}, \qquad
\pg_{\gamma}=\pg\cap\ggo_{\gamma}.
$$
If $p:\pg\mapsto\pg_{\gamma}$ is the orthogonal projection relative
to $\ip_{\pg}$, then it is easy to see that
$$
m_{\gamma}:V\mapsto\pg_{\gamma}, \qquad m_{\gamma}=p\circ m,
$$
is precisely the moment map for the action of $G_{\gamma}$ on $V$.
}
\end{example}

In the cases considered in detail in this paper we will have
$(G_{\gamma},K_{\gamma})$ equal to
$(\Spe(\frac{n}{2},\RR),\U(\frac{n}{2}))$ (symplectic),
$(\Gl(\frac{n}{2},\CC),\U(\frac{n}{2}))$ (complex),
$(\Gl(\frac{n}{4},\HH),\Spe(\frac{n}{4}))$ (hypercomplex) and
$(\Gl(n,\RR),\Or(n))$ ($\gamma=0$).

\section{Variety of compatible metrics}\label{var}

Let us consider as a parameter space for the set of all real
nilpotent Lie algebras of a given dimension $n$, the set $\nca$ of
all nilpotent Lie brackets on a fixed $n$-dimensional real vector
space $\ngo$.  If
$$
V=\lam=\{\mu:\ngo\times\ngo\mapsto\ngo : \mu\; \mbox{skew-symmetric
bilinear map}\},
$$
then
$$
\nca=\{\mu\in V:\mu\;\mbox{satisfies Jacobi and is nilpotent}\}
$$
is an algebraic subset of $V$.  Indeed, the Jacobi identity and the
nilpotency condition are both determined by zeroes of polynomials.

We fix a tensor $\gamma$ on $\ngo$ (or a set of tensors), and let
$G_{\gamma}$ denote the subgroup of $\Gl(n)$ preserving $\gamma$.
These groups act naturally on $V$ by (\ref{action}) and leave $\nca$
invariant.  Consider the subset $\nca_{\gamma}\subset\nca$ given by
$$
\nca_{\gamma}=\{\mu\in\nca:{\rm IC}(\gamma,\mu)=0\},
$$
that is, those nilpotent Lie brackets for which $\gamma$ is
integrable (see (\ref{closedg})). $\nca_{\gamma}$ is also an
algebraic variety since ${\rm IC}(\gamma,\mu)$ is always linear on
$\mu$.  Recall that
$$
W_{\gamma}=\{\mu\in V:{\rm IC}(\gamma,\mu)=0\}
$$
is a $G_{\gamma}$-invariant linear subspace of $V$, and
$\nca_{\gamma}=\nca\cap W_{\gamma}$.

For each $\mu\in\nca$, let $N_{\mu}$ denote the simply connected
nilpotent Lie group with Lie algebra $(\ngo,\mu)$.  We now consider
an identification of each point of $\nca_{\gamma}$ with a compatible
metric.  Fix an inner product $\ip$ on $\ngo$ compatible with
$\gamma$, that is, such that (\ref{ortcon}) holds. We identify each
$\mu\in\nca_{\gamma}$ with a class-$\gamma$ metric structure
\begin{equation}\label{ideg}
\mu\longleftrightarrow   (N_{\mu},\gamma,\ip),
\end{equation}
where all the structures are defined by left invariant translation.
Therefore, each $\mu\in\nca_{\gamma}$ can be viewed in this way as a
metric compatible with the class-$\gamma$ nilpotent Lie group
$(N_{\mu},\gamma)$, and two metrics $\mu,\lambda$ are compatible
with the same geometric structure if and only if they live in the
same $G_{\gamma}$-orbit. Indeed, the action of $G_{\gamma}$ on
$\nca_{\gamma}$ has the following interpretation: each $\vp\in
G_{\gamma}$ determines a Riemannian isometry preserving the
geometric structure
$$
(N_{\vp.\mu},\gamma,\ip)\mapsto
(N_{\mu},\gamma,\la\vp\cdot,\vp\cdot\ra)
$$
by exponentiating the Lie algebra isomorphism
$\vp^{-1}:(\ngo,\vp.\mu)\mapsto(\ngo,\mu)$.  We then have the
identification $G_{\gamma}.\mu=\cca(N_{\mu},\gamma)$, and more in
general the following

\begin{proposition}\label{upto}
Every class-$\gamma$ metric structure $(N',\gamma',\ip')$ on a
nilpotent Lie group $N'$ of dimension $n$ is isometric-isomorphic to
a $\mu\in\nca_{\gamma}$.
\end{proposition}

\begin{proof}
We can assume that the Lie algebra of $N'$ is $(\ngo,\lambda)$ for
some $\lambda\in\nca$.  There exist $\vp\in\Gl(n)$ and
$\psi\in\Or(\ngo,\ip)$ such that $\vp.\ip'=\ip$ and
$\psi(\vp.\gamma')=\gamma$.  Thus the Lie algebra isomorphism
$\psi\vp:(\ngo,\lambda)\mapsto(\ngo,\mu)$, where
$\mu=\psi\vp.\lambda$, satisfies $\psi\vp.\ip'=\ip$ and
$\psi\vp.\gamma'=\gamma$ and so it defines an isometry
$$
(N',\gamma',\ip')\mapsto(N_{\mu},\gamma,\ip)
$$
which is also an isomorphism between the class-$\gamma$ nilpotent
Lie groups $(N',\gamma')$ and $(N,\gamma)$, concluding the proof.
\end{proof}

According to the above proposition and identification (\ref{ideg}),
the orbit $G_{\gamma}.\mu$ parameterizes all the left invariant
metrics which are compatible with $(N_{\mu},\gamma)$ and hence we
may view $\nca_{\gamma}$ as the space of all class-$\gamma$ metric
structures on nilpotent Lie groups of dimension $n$. Since two
metrics $\mu,\lambda\in\nca_{\gamma}$ are isometric if and only if
they live in the same $K_{\gamma}$-orbit, where
$K_{\gamma}=G_{\gamma}\cap\Or(\ngo,\ip)$ (see Appendix \ref{app}),
we have that $\nca_{\gamma}/K_{\gamma}$ parameterizes class-$\gamma$
metric nilpotent Lie groups of dimension $n$ up to isometry and
$G_{\gamma}.\mu/K_{\gamma}$ do the same for all the compatible
metrics on $(N_{\mu},\gamma)$.

We now recall a crucial fact which is the link between the study of
left invariant compatible metrics for geometric structures on
nilpotent Lie groups and the results from real geometric invariant
theory exposed in Section \ref{git}.  The interplay is based on the
identification given in (\ref{ideg}), and it will be used in the
proofs of the remaining results of this section.  The proof of the
following proposition follows just from a simple comparison between
formulas (\ref{ricci}) and (\ref{defm2}).

\begin{proposition}\label{momric}
Let $m:V\mapsto\pg$ and $m_{\gamma}:V\mapsto\pg_{\gamma}$  be the
moment maps for the actions of $\Gl(n)$ and $G_{\gamma}$ on
$V=\lam$, respectively {\rm (see Examples \ref{gln}, \ref{momg})},
where $\pg$ is the space of symmetric maps of $(\ngo,\ip)$ and
$\pg_{\gamma}$ the subspace of those leaving $\gamma$ invariant.
\begin{itemize}
\item[(i)] For each $\mu\in\nca\subset V$,
$$
m(\mu)=8\Ric_{\mu},
$$
where $\Ric_{\mu}$ is the Ricci operator of the Riemannian manifold
$(N_{\mu},\ip)$.

\item[(ii)] For each $\mu\in\nca_{\gamma}\subset V$,
$$
m_{\gamma}(\mu)=8\Ricg_{\mu},
$$
where $\Ricg_{\mu}$ is the invariant Ricci operator of
$(N_{\mu},\gamma,\ip)$, that is, the orthogonal projection of the
Ricci operator $\Ric_{\mu}$ on $\pg_{\gamma}$.
\end{itemize}
\end{proposition}

Let us now go back to our search for the best compatible metric.
The identification (\ref{ideg}) allows us to view each point of the
variety $\nca_{\gamma}$ as a class-$\gamma$ metric structure on a
nilpotent Lie group of dimension $n$. In this light, it is natural
to consider the functional $F:\nca_{\gamma}\mapsto\RR$ given by
$F(\mu)=\tr(\Ricg_{\mu})^2$, which in some sense measures how far is
the metric $\mu$ from having $\Ricg_{\mu}=0$, which is the goal
proposed in Section \ref{geometric} (see (\ref{uncondition})).  The
critical points of $F$ may be therefore considered compatible
metrics of particular significance.  However, we should consider
some normalization since $F(t\mu)=t^4F(\mu)$ for all $t\in\RR$.

For any $\mu\in\nca_{\gamma}$ we have that $\scalar(\mu)=-\unc
||\mu||^2$ (see (\ref{ricci})), which says that normalizing by
scalar curvature and by the spheres of $V$ is equivalent:
$$
\{\mu\in\nca_{\gamma}:\scalar(\mu)=s\}=\{\mu\in\nca_{\gamma}:||\mu||^2=-4s\},
\quad \forall s<0.
$$
The critical points of $F:\PP V\mapsto\RR$,
$F([\mu])=\tr(\Ricg_{\mu})^2/||\mu||^4$, which lie in
$\PP\nca_{\gamma}=\pi(\nca_{\gamma})$ appears then as very natural
candidates, since it is like we are restricting $F$ to the subset of
all class-$\gamma$ metric structures having a given scalar
curvature.

It follows from Proposition \ref{momric}, (ii), that
$$
F([\mu])=\frac{1}{64}||m_{\gamma}([\mu])||^2,
$$
where $m_{\gamma}:\PP V\mapsto\pg_{\gamma}$ is the moment map for
the action of $G_{\gamma}$ on $\PP V$.  We then obtain from Lemma
\ref{marian2} and (\ref{actiong}) that
$$
\grad(F)_{[\mu]}=-\frac{1}{16}\pi^*\delta_{\mu}(\Ricg_{\mu}), \qquad
||\mu||=1,
$$
where $\pi^*:V\mapsto\tang_{[\mu]}\PP V$ denotes the derivative of
the projection map $\pi:V\mapsto\PP V$.  This shows that
$[\mu]\in\PP V$ is a critical point of $F$ if and only if
$\Ricg_{\mu}=cI+D$ for some $c\in\RR$ and $D\in\Der(\mu)$
($=\Ker{\delta_{\mu}}$).  By applying Theorem \ref{marian} to our
situation we obtain the main result of this paper.

\begin{theorem}\label{equiv1g}
Let $F:\PP V\mapsto\RR$ be defined by
$F([\mu])=\tr(\Ricg_{\mu})^2/||\mu||^4$.  Then for $\mu\in V$ the
following conditions are equivalent:
\begin{itemize}
\item[(i)] $[\mu]$ is a critical point of $F$.

\item[(ii)] $F|_{G_{\gamma}.[\mu]}$ attains its minimum value at $[\mu]$.

\item[(iii)] $\Ricg_{\mu}=cI+D$ for some $c\in\RR$, $D\in\Der(\mu)$.
\end{itemize}
Moreover, all the other critical points of $F$ in the orbit
$G_{\gamma}.[\mu]$ lie in $K_{\gamma}.[\mu]$.
\end{theorem}

We now rewrite the above result in geometric terms, by using the
identification (\ref{ideg}), Proposition \ref{soliton1g} and
Definition \ref{minimalg}.

\begin{theorem}\label{equiv2g}
Let $(N,\gamma)$ be a nilpotent Lie group endowed with an invariant
geometric structure $\gamma$.  Then the following conditions on a
left invariant Riemannian metric $\ip$ which is compatible with
$(N,\gamma)$ are equivalent:
\begin{itemize}
\item[(i)] $\ip$ is minimal.

\item[(ii)] $\ip$ is an invariant Ricci soliton.

\item[(iii)] $\Ricg_{\ip}=cI+D$ for some $c\in\RR$, $D\in\Der(\ngo)$.
\end{itemize}
Moreover, there is at most one compatible left invariant metric on
$(N,\gamma)$ up to isometry  (and scaling) satisfying any of the
above conditions.
\end{theorem}

Recall that the proof of this theorem does not use the integrability
of $\gamma$, and so it is valid for the `almost' versions as well.

We also note that part (i) of Theorem \ref{equiv1g} makes possibly
the study of minimal compatible metrics by a variational method.
Indeed, the projective algebraic variety $\PP\nca_{\gamma}$ may be
viewed as the space of all class-$\gamma$ metric structures on
$n$-dimensional nilpotent Lie groups with a given scalar curvature,
and those which are minimal are precisely the critical points of
$F:\PP\nca_{\gamma}\mapsto\RR$.  This variational approach will be
used quite often in the search for explicit examples in the
following section and in \cite{praga}.

The theorems above propose then as privileged these compatible
metrics called minimal, which have a neat characterization (see
(iii)), are critical points of a natural curvature functional
(square norm of Ricci), minimize such a functional when restricted
to the compatible metrics for a given geometric structure, and are
solitons for a natural evolution flow.  Moreover, the uniqueness up
to isometry of such special metrics holds. But a weakness of this
approach is the existence problem.  How special are the symplectic
or (almost-) complex structures admitting a minimal metric?.  So
far, we know how to deal with this `existence question' only by
giving several examples, which is the goal of Section \ref{exa} and
\cite{praga}.  The only obstruction we have found is in the case
$\RR I\not\subset\ggo_{\gamma}$, namely when $\Der(\ngo)$ is
nilpotent.  These Lie algebras are called {\it characteristically
nilpotent} and have been extensively studied in the last years.
Examples of characteristically nilpotent Lie algebra admitting a
symplectic structure have been recently given in \cite{Brd}, and
these are the only counterexamples we know to the existence problem.

\begin{corollary}\label{isoiso}
Let $\gamma,\gamma'$ be two geometric structures on a nilpotent Lie
group $N$, and assume that they admit minimal compatible metrics
$\ip$ and $\ip'$, respectively.  Then $\gamma$ is isomorphic to
$\gamma'$ if and only if there exists $\vp\in\Aut(\ngo)$ and $c>0$
such that $\gamma'=\vp.\gamma$ and
$$
\la\vp X,\vp Y\ra'=c\la X,Y\ra \qquad \forall X,Y\in\ngo.
$$
In particular, if $\gamma$ and $\gamma'$ are isomorphic then their
respective minimal compatible metrics are necessarily isometric up
to scaling (recall that $c=1$ when $\scalar(\ip)=\scalar(\ip')$).
\end{corollary}

We have here a very useful tool to distinguish two geometric
structures.  Indeed, the corollary allows us to do it by looking at
their respective minimal compatible metrics, that is, with
Riemannian data.  This is a remarkable advantage since we suddenly
have a great deal of invariants.  This method will be used in the
next section to find explicit continuous families of pairwise
non-isomorphic geometric structures in low dimensions, mainly by
using only one Riemannian invariant: the eigenvalues of the Ricci
operator.

\begin{remark}\label{liegroups}
{\rm The Ricci curvature operator of a left invariant metric $\ip$
on a Lie group $G$ is given by
$$
\Ric_{\ip}=R_{\ip}-\unm B_{\ip}-D_{\ip},
$$
where $R_{\ip}$ is defined by (\ref{ricci}), $B_{\ip}$ is the
Killing form of the Lie algebra $\ggo$ of $G$ in terms of $\ip$,
$D_{\ip}$ is the symmetric part of $\ad{Z_{\ip}}$ and
$Z_{\ip}\in\ggo$ is defined by $\la Z_{\ip},X\ra=\tr(\ad{X})$ for
any $X\in\ggo$.  Recall that $Z_{\ip}=0$ if and only if $\ggo$ is
unimodular, and that $\Ric_{\ip}=R_{\ip}$ in the nilpotent case.  If
we consider the tensor $R$ instead of the Ricci tensor, and the
variety of all Lie algebras $\lca$ rather than just the nilpotent
ones, to define and state all the notions, flows, identifications
and results in Sections \ref{geometric} and \ref{var}, then
everything is still valid for Lie groups in general, with the only
exception of the first part of Corollary \ref{isoiso}. We only have
to consider the corresponding invariant part $R^{\gamma}$ and
replace $\scalar(\ip)$ with $\tr{R_{\ip}}$ each time it appears.
The only detail to be careful with is that if two
$\mu,\lambda\in\lca_{\gamma}$ lie in the same $K_{\gamma}$-orbit
then they are isometric, but the converse might not be true.  Recall
that the uniqueness result in Theorem \ref{equiv2g} is nevertheless
valid.

The reason why we decided to work only in the nilpotent case is
that, at least at first sight, the use of this `unnamed' tensor $R$
make minimal and soliton metrics, as well as the functionals and
evolution flows, into concepts lacking in geometric sense.  For
instance, we have found ourselves with the unpleasant fact that some
K$\ddot{{\rm a}}$hler metrics on solvable Lie groups would not be
minimal viewed as compatible metrics for the corresponding
symplectic structures, in spite of $\Ricac_{\ip}=0$.

For a compact simple Lie group, $-B$ is minimal for the case
$\gamma=0$, and if $\ggo=\kg\oplus\pg$ is the Cartan decomposition
for a non-compact semi-simple Lie algebra $\ggo$, then it is easy to
see that the metric $\ip$ given by $\la\kg,\pg\ra=0$,
$\ip|_{\kg\times\kg}=-B$ and $\ip|_{\pg\times\pg}=B$ is minimal as
well.    }
\end{remark}

\section{Applications}\label{exa}

\subsection{Symplectic structures}\label{nilpsymp}
Let $N$ be a real $2n$-dimensional nilpotent Lie group with Lie
algebra $\ngo$, whose Lie bracket is denoted by $\mu
:\ngo\times\ngo\mapsto\ngo$.  An invariant {\it symplectic}
structure on $N$ is defined by a $2$-form $\omega$ on $\ngo$
satisfying
$$
\omega(X,\cdot)\equiv 0 \quad \mbox{if and only if} \quad X=0 \quad
(\mbox{non-degenerate}),
$$
and for all $X,Y,Z\in\ngo$
\begin{equation}\label{closed}
\omega(\mu(X,Y),Z)+\omega(\mu(Y,Z),X)+\omega(\mu(Z,X),Y)=0 \quad
({\rm closed}, \dif\omega=0).
\end{equation}
Fix a symplectic nilpotent Lie group $(N,\omega)$.  A left invariant
Riemannian metric which is compatible with $(N,\omega)$ is
determined by an inner product $\ip$ on $\ngo$ such that if
$$
\omega(X,Y)=\la X,J_{\ip}Y\ra\quad\forall\; X,Y\in\ngo\quad{\rm
then}\quad J_{\ip}^2=-I.
$$
For the geometric structure $\gamma=\omega$ we have that
$$
G_{\gamma}=\Spe(n,\RR)=\{ g\in\Gl(2n):g^tJg=J\}, \qquad
K_{\gamma}=\U(n),
$$
and the Cartan decomposition of $\ggo_{\gamma}=\spg(n,\RR)=\{
A\in\glg(2n):A^tJ+JA=0\}$ is given by
$$
\spg(n,\RR)=\ug(n)\oplus\pg_{\gamma}, \qquad \pg_{\gamma}=\{
A\in\pg:AJ=-JA\}.
$$
Thus the invariant Ricci tensor $\riccig$ coincides with the
anti-complexified Ricci tensor (see \cite{LeWng}) and for any
$\ip\in\cca$,
\begin{equation}\label{ac-cop}
\Ricg_{\ip}=\Ricac_{\ip}=\unm(\Ric_{\ip}+J_{\ip}\Ric_{\ip}J_{\ip}).
\end{equation}
This implies that our `goal' condition $\Ricg_{\ip}=0$ (see
(\ref{uncondition})) is equivalent to have hermitian Ricci tensor.
Also, the evolution flow considered in Section \ref{geometric} is
not other than the anti-complexified Ricci flow.

Concerning the search for the best compatible left invariant metric
for a symplectic nilpotent Lie group, our first result is negative.

\begin{proposition}\label{noherm}
Let $(N,\omega)$ be a symplectic nilpotent Lie group.  Then
$(N,\omega)$ does not admit any compatible left invariant metric
with hermitian Ricci tensor, unless $N$ is abelian.
\end{proposition}

\begin{proof}
We first note that since $\mu$ is nilpotent the center $\zg$ of
$(\ngo,\mu)$ meets non-trivially the derived Lie subalgebra
$\mu(\ngo,\ngo)$, unless $\mu=0$ (i.e. $\ngo$ abelian).  Assume that
$\ip\in\cca(N,\omega)$ has hermitian Ricci tensor and consider the
orthogonal decomposition $\ngo=\vg\oplus\mu(\ngo,\ngo)$.  If
$Z\in\zg$ then $JZ\in\vg$. In fact, it follows from (\ref{closed})
that
$$
\la\mu(X,Y),JZ\ra=\omega(\mu(X,Y),Z)=0 \qquad \forall\; X,Y\in\ngo.
$$
Now, the above equation, the fact that $\Ric_{\ip}J=J\Ric_{\ip}$ and
the definition of $\Ric_{\ip}$ (see (\ref{ricci})) imply that
$$
0\leq\la\Ric_{\ip}Z,Z\ra=\la\Ric_{\ip}JZ,JZ\ra\leq 0,
$$
and hence
$$
\unc\sum_{ij}\la\mu(X_i,X_j),Z\ra^2=\la\Ric_{\ip}Z,Z\ra=0 \qquad
\forall\; Z\in\zg.
$$
Thus $\mu(\ngo,\ngo)\perp\zg$ and so $\ngo$ must be abelian by the
observation made in the beginning of the proof.
\end{proof}

We now review the variational approach developed in Section
\ref{var}.  Fix a non-degenerate $2$-form $\omega$ on $\ngo$, and
let $\Spe(n,\RR)$ denote the subgroup of $\Gl(2n)$ preserving
$\omega$, that is,
$$
\Spe(n,\RR)=\{ \vp\in\Gl(2n):\omega(\vp X,\vp Y)=\omega(X,Y)
\quad\forall\; X,Y\in\ngo\}.
$$
Consider the algebraic subvariety $\nca_s:=\nca_{\gamma}\subset\nca$
given by
$$
\nca_s=\{\mu\in\nca:\dif_{\mu}\omega=0\},
$$
that is, those nilpotent Lie brackets for which $\omega$ is closed
(see (\ref{closed})).  By fixing an inner product $\ip$ on $\ngo$
satisfying that
$$
\omega=\la\cdot,J\cdot\ra\qquad{\rm with}\qquad J^2=-I,
$$
(\ref{ideg}) identify each $\mu\in\nca_s$ with the
almost-K$\ddot{{\rm a}}$hler manifold $(N_{\mu},\omega,\ip,J)$.  The
action of $\Spe(n,\RR)$ on $\nca_s$ has the following
interpretation:  each $\vp\in\Spe(n,\RR)$ determines a Riemannian
isometry which is also a symplectomorphism
$$
(N_{\vp.\mu},\omega,\ip,J)\mapsto
(N_{\mu},\omega,\la\vp\cdot,\vp\cdot\ra,\vp^{-1}J\vp)
$$
by exponentiating the Lie algebra isomorphism
$\vp^{-1}:(\ngo,\vp.\mu)\mapsto(\ngo,\mu)$.

Let $\ngo$ be a $6$-dimensional vector space with basis $\{
X_1,...,X_{6}\}$ over $\RR$, and consider the non-degenerate
$2$-form
$$
\omega=\alpha_1\wedge\alpha_{6}+\alpha_2\wedge\alpha_{5}+\alpha_3\wedge\alpha_{4},
$$
where $\{\alpha_1,...,\alpha_{6}\}$ is the dual basis of $\{ X_i\}$.
For the compatible inner product $\la X_i,X_j\ra=\delta_{ij}$ we
have that $\omega=\la\cdot,J\cdot\ra$ for
$$
J=\left[\begin{smallmatrix} &&&&&-1\\ &0&&&\cdot&\\ &&&-1&&\\ &&1&&&\\ &\cdot&&&0&\\
1&&&&&\end{smallmatrix}\right].
$$
In the following example the symplectic structure will always be
$\omega$, the almost-complex structure $J$ and the compatible metric
$\ip$.  We will vary Lie brackets and use constantly identification
(\ref{ideg}).

\begin{example}\label{m26}
{\rm We shall take advantage of the variational nature of Theorem
\ref{equiv1g} to find explicit examples of minimal compatible
metrics.  Consider for each $6$-upla $\{ a,...,f\}$ of real numbers
the skew-symmetric bilinear form $\mu=\mu(a,b,c,d,e,f)\in
V=\Lambda^2(\RR^6)^*\otimes\RR^6$ defined by
$$
\begin{array}{l}
\mu(X_1,X_2)=aX_3, \qquad\mu(X_1,X_3)=bX_4, \qquad \mu(X_1,X_4)=cX_5,\\

\mu(X_1,X_5)=dX_6, \qquad\mu(X_2,X_3)=eX_5,
\qquad\mu(X_2,X_4)=fX_6.
\end{array}
$$

Our plan is to find first the critical points of $F$ restricted to
the set $\{[\mu(a,...,f)]:a,...,f\in\RR\}$ and after that to show by
using the characterization given in part (iii) of the theorem that
they are really critical points of $F:\PP V\mapsto\RR$.  We can see
by a simple computation that $\Ricac_{\mu}$ is given by the diagonal
matrix with entries
$$
\begin{array}{rl}
\Ricac_{\mu}&=-\unc\left[\begin{smallmatrix} a^2+b^2+c^2+2d^2+f^2\\ a^2+c^2-d^2+2e^2+f^2\\
-a^2+2b^2-c^2+e^2-f^2\\
a^2-2b^2+c^2-e^2+f^2\\ -a^2-c^2+d^2-2e^2-f^2\\
-a^2-b^2-c^2-2d^2-f^2\end{smallmatrix}\right],
\end{array}
$$
and hence we are interested in the critical points of
$$
\begin{array}{rl}
F(\mu)=&\tr(\Ricac_{\mu})^2=F(a,...,f)\\
=&\frac{1}{8}\Big((a^2+b^2+c^2+2d^2+f^2)^2+(a^2+c^2-d^2+2e^2+f^2)^2 \\
&+(-a^2+2b^2-c^2+e^2-f^2)^2\Big)
\end{array}
$$
restricted to any leaf of the form $a^2+b^2+c^2+d^2+e^2+f^2\equiv$
const., which are easily seen to depend of three parameters.  We
still have to impose the Jacobi and closeness conditions on these
critical points (or equivalently to find the intersection with
$\nca_s$), after which we obtain the following ellipse of symplectic
structures:
$$
\{\mu_{xy}=\mu(x,1,x+y,1,1,y):x^2+y^2+xy=1\}.
$$
It follows from the formula for $\Ricac_{\mu}$ given above that
$$
\Ricac_{\mu_{xy}}=-\unc\left[\begin{smallmatrix} 5&&&&&\\ &3&&&&\\ &&1&&&\\
&&&-1&&\\ &&&&-3&\\ &&&&&-5\end{smallmatrix}\right]
=-\frac{7}{4}I+\unm\left[\begin{smallmatrix} 1&&&&&\\ &2&&&&\\ &&3&&&\\
&&&4&&\\ &&&&5&\\ &&&&&6 \end{smallmatrix}\right]\in\RR
I+\Der(\mu_{xy}),
$$
showing definitely that this is a curve of minimal compatible
metrics.  We furthermore have that the Ricci tensor of the metrics
$\mu_{xy}$ is given by
$$
\Ric_{\mu_{xy}}=-\unc\left[\begin{smallmatrix} 4-y^2&&&&&\\ &2-xy&&&&\\ &&2-x^2&&&\\ &&&1-x^2&&\\ &&&&-1-xy&\\
&&&&&-1-y^2\end{smallmatrix}\right],
$$
which clearly shows that they are pairwise non-isometric for
$x,y\geq 0$.  It then follows from Corollary \ref{isoiso} that
$$
\{(N_{\mu_{xy}},\omega):x^2+y^2+xy=1,\; x,y\geq 0\}
$$
is a curve of pairwise non-isomorphic symplectic nilpotent Lie
groups.  There are three $6$-dimensional nilpotent Lie groups
involved, $N_{\mu_{10}}$, $N_{\mu_{01}}$ and $N_{\mu_{xy}}$,
$x,y>0$, denoted in \cite[Table A.1]{Slm} by $(0,0,12,13,14,23+15)$,
$(0,0,0,12,14-23,15+34)$ and $(0,0,12,13,14+23,24+15)$,
respectively.  }
\end{example}

\subsection{Complex structures}\label{complex}
Let $N$ be a real $2n$-dimensional nilpotent Lie group with Lie
algebra $\ngo$, whose Lie bracket is denoted by $\mu
:\ngo\times\ngo\mapsto\ngo$.  An invariant {\it almost-complex}
structure on $N$ is defined by a map $J:\ngo\mapsto\ngo$ satisfying
$J^2=-I$.  If in addition $J$ satisfies the integrability condition
\begin{equation}\label{integral}
\mu(JX,JY)=\mu(X,Y)+J\mu(JX,Y)+J\mu(X,JY), \qquad \forall
X,Y\in\ngo,
\end{equation}
then $J$ is said to be a {\it complex} structure.  A left invariant
Riemannian metric which is {\it compatible} with $(N,J)$, also
called an {\it almost-hermitian metric}, is given by an inner
product $\ip$ on $\ngo$ such that
$$
\la JX,JY\ra=\la X,Y\ra \qquad \forall X,Y\in\ngo.
$$
As in the symplectic  case, condition $\Ricc_{\ip}=0$ is forbidden
for non-abelian $N$ since $\tr{\Ricc_{\ip}}=\scalar(\ip)<0$.

Let $\{ X_1,...,X_4,Z_1,Z_2\}$ be a basis for $\ngo$.  The complex
structure and the compatible metric in the following example will
always be defined by
$$
\begin{array}{lll}
J=\left[\begin{smallmatrix} 0&-1&&&&\\ 1&0&&&&\\ &&0&-1&&\\
&&1&0&&\\ &&&&0&-1\\ &&&&1&0
\end{smallmatrix}\right], && \la X_i,X_j\ra=\la Z_i,Z_j\ra=\delta_{ij}.
\end{array}
$$

\begin{example} {\rm The following curve has been obtained via the
variational method provided by Theorem \ref{equiv1g}, by using an
approach very similar to that in Example \ref{m26}.  Let $\mu_t$ be
the curve defined by {\small
$$
\begin{array}{lll}
\mu(X_1,X_3)=-tsZ_2, && \mu(X_2,X_3)=sZ_1, \\
 \mu(X_1,X_4)=sZ_1, && \mu(X_2,X_4)=s(2-t)Z_2,
\end{array} \qquad s=\sqrt{2+t^2+(2-t)^2}, \quad t\in\RR.
$$}
It is easy to check that $(N_{\mu_t},J)$ is a non-abelian complex
nilpotent Lie group for all $t\in\RR$. Moreover,
$\Ric_{\mu}|_{\ngo_1}$ is diagonal and hence both
$\Ricc_{\mu}|_{\ngo_1}$ and $\Ricc_{\mu}|_{\ngo_2}$ are scalar
multiples of the identity (recall that the invariant Ricci operator
is given in this case by
$\Ricg_{\ip}=\Ricc_{\ip}=\unm(\Ric_{\ip}-J\Ric_{\ip}J)$).  Since
$\mu$ is two-step nilpotent, this easily implies that $\ip$ is a
minimal compatible metric for all $(N_{\mu_t},J)$.   Indeed, if
$\Ricc_{\mu}|_{\ngo_1}=pI$ and $\Ricc_{\mu}|_{\ngo_2}=qI$ for some
$p,q\in\RR$, we would have that
\begin{equation}\label{multide}
\Ricc_{\mu}=\left[\begin{smallmatrix}pI&\\&qI\end{smallmatrix}\right]=(2p-q)I+
\left[\begin{smallmatrix}(q-p)I&\\
&2(q-p)I\end{smallmatrix}\right]\in\RR I+\Der(\mu).
\end{equation}
It follows from
$$
\Ric_{\mu_{t}}|_{\ngo_2}=\unm\left[\begin{smallmatrix} 2s^2&0\\
0&s^2(t^2+(2-t)^2)\end{smallmatrix}\right],
$$
that the hermitian manifolds $\{(N_{\mu_t},J,\ip):1\leq t<\infty\}$
are pairwise non-isometric since
$$
s^2(t^2+(2-t)^2)-2s^2=(t^2+(2-t)^2)^2-4
$$
is a strictly increasing non-negative function for $1\leq t$, which
vanishes if and only if $t=1$.  We therefore obtain a curve
$\{(N_{\mu_t},J):1\leq t<\infty\}$ of pairwise non-isomorphic
non-abelian complex nilpotent Lie groups.  A natural question is
which are the nilpotent Lie groups involved.  We have for all $t$
that
$$
j_{\mu_t}(Z_1)=\left[\begin{smallmatrix} &&0&-s\\ &&-s&0\\ 0&s&&\\
s&0&&\end{smallmatrix}\right],\qquad
j_{\mu_t}(Z_2)=\left[\begin{smallmatrix} &&ts&0\\ &&0&-(2-t)s\\
-ts&0&&\\ 0&(2-t)s&&\end{smallmatrix}\right],
$$
(see Appendix \ref{app}) and hence $j_{\mu_t}(Z)$ is non-singular if
and only if
$$
-t(2-t)\la Z,Z_1\ra^2-\la Z,Z_2\ra^2\ne 0.
$$
This implies that $\mu_t$ is isomorphic to the complex Heisenberg
Lie algebra (i.e. when $j_{\mu_t}(Z)$ is non-singular for any
non-zero $Z\in\ngo_2$) if and only if $1\leq t<2$, providing a curve
on the Iwasawa manifold.  Furthermore, $(N_{\mu_1},J)$ is the
bi-invariant complex structure and it can be showed by computing
$j_{\mu_t}(Z)^2$ that $(N_{\mu_t},\ip)$ is not modified H-type for
any $1<t$.  We finally note that $\mu_2$ is isomorphic to the group
denoted by $(0,0,0,0,12,14+23)$ in \cite{Slm}, and one can easily
see by discarding any other possibility that actually
$\mu_t\simeq\mu_2$ for all $2\leq t<\infty$, which gives rise a
curve of pairwise non-isomorphic structures on such a group.  }
\end{example}

\subsection{Hypercomplex structures}\label{hyper}
Let $N$ be a real $4n$-dimensional nilpotent Lie group with Lie
algebra $\ngo$, whose Lie bracket is denoted by $\mu
:\ngo\times\ngo\mapsto\ngo$.  An invariant {\it hypercomplex}
structure on $N$ is defined by a triple $\{ J_1,J_2,J_3\}$ of
complex structures on $\ngo$ (see Section \ref{complex}) satisfying
the quaternion identities
\begin{equation}\label{quat}
J_i^2=-I,\quad i=1,2,3, \qquad J_1J_2=J_3=-J_2J_1.
\end{equation}
An inner product $\ip$ on $\ngo$ is said to be {\it compatible} with
$\{ J_1,J_2,J_3\}$, also called an {\it hyper-hermitian metric}, if
\begin{equation}\label{ortconhyp}
\la J_iX,J_iY\ra=\la X,Y\ra \qquad \forall X,Y\in\ngo, \; i=1,2,3.
\end{equation}
Two hypercomplex nilpotent Lie groups $(N,\{ J_1,J_2,J_3\})$ and
$(N',\{ J_1',J_2',J_3'\})$ are said to be {\it isomorphic} if there
exists an isomorphism $\vp:\ngo'\mapsto\ngo$ such that
$$
\vp J_i'\vp^{-1}=J_i, \quad i=1,2,3.
$$
There are no non-abelian nilpotent Lie groups of dimension $4$
admitting an hypercomplex structure.  In dimension $8$, hypercomplex
nilpotent Lie groups have been determined by I. Dotti and A. Fino in
\cite{DttFin0} and \cite{DttFin2}.  They proved the following strong
restrictions on an $8$-dimensional nilpotent Lie algebra $\ngo$
which admits an hypercomplex structure: $\ngo$ has to be $2$-step
nilpotent, $\dim{\mu(\ngo,\ngo)}\leq 4$, there exists a
decomposition $\ngo=\ngo_1\oplus\ngo_2$ such that $\dim{\ngo_i}=4$,
$\ngo_i$ is $\{ J_1,J_2,J_3\}$-invariant and
$\mu(\ngo,\ngo)\subset\ngo_2\subset\zg$, where $\zg$ is the center
of $\ngo$.  Thus the Lie bracket of $\ngo$ is just given by a
skew-symmetric bilinear form $\mu:\ngo_1\times\ngo_1\mapsto\ngo_2$,
and those for which a fixed $\{ J_1,J_2,J_3\}$ is integrable are
also completely described in \cite{DttFin2} (see also
\cite{DttFin1}) as a $16$-dimensional subspace $W_h$ of the
$24$-dimensional vector space $W=\Lambda^2\ngo_1^*\otimes\ngo_2$ of
all such forms.  If we ask in addition abelian (i.e.
$\mu=\mu(J_i\cdot,J_i\cdot)$), then we get a subspace $W_{ah}$ of
dimension $12$.

What shall be studied here are the isomorphism classes of such
structures and the existence of minimal compatible metrics.

Fix basis $\{ X_1,X_2,X_3,X_4\}$ and $\{ Z_1,Z_2,Z_3,Z_4\}$ of
$\ngo_1$ and $\ngo_2$, respectively.  The compatible metric will be
$\la X_i,X_j\ra=\la Z_i,Z_j\ra=\delta_{ij}$ and the hypercomplex
structure will always act on $\ngo_i$ by
$$
J_1=\left[\begin{smallmatrix} 0&-1&&\\ 1&0&&\\ &&0&-1\\
&&1&0\end{smallmatrix}\right], \quad J_2=\left[\begin{smallmatrix}
&&-1&0\\ &&0&1\\ 1&0&&\\ 0&-1&&\end{smallmatrix}\right],\quad
J_3=\left[\begin{smallmatrix} &&&-1\\ &&-1&\\ &1&&\\
1&&&\end{smallmatrix}\right].
$$

\begin{proposition}
{\rm (i)} Every hypercomplex $8$-dimensional nilpotent Lie group
admits a minimal compatible metric.

\no {\rm (ii)} Two hypercomplex $8$-dimensional nilpotent Lie groups
$(N_{\mu},\{ J_1,J_2,J_3\})$ and $(N_{\lambda},\{ J_1,J_2,J_3\})$
are isomorphic if and only if $\mu$ and $c\lambda$ lie in the same
$\Spe(1)\times\Spe(1)$-orbit for some non-zero $c\in\RR$.

\no {\rm (iii)} The moduli space of all $8$-dimensional hypercomplex
nilpotent Lie groups up to isomorphism is parameterized by
$$
\PP W_{h}/\Spe(1)\times\Spe(1).
$$
The representation $W_h$ is equivalent to
$(\sug(2)\otimes\RR^4)\oplus\RR^4$, where $\sug(2)$ is the adjoint
representation and $\RR^4$ is the standard representation of
$\SU(2)=\Spe(1)$ viewed as real.  Since the isotropy of an element
in general position is finite, the dimension of this quotient is
$15-6=9$.

\no {\rm (iv)} The moduli space of all $8$-dimensional abelian
hypercomplex nilpotent Lie groups up to isomorphism is parameterized
by
$$
\PP W_{ah}/\Spe(1)\times\Spe(1).
$$
The representation $W_{ah}$ is equivalent to $\sug(2)\otimes\RR^4$,
and since the isotropy of an element in general position is again
finite, the dimension of this quotient is $11-6=5$.
\end{proposition}

\begin{proof}
Since the only symmetric transformations of $\ngo_i=\RR^4$ commuting
with all the $J_i's$ are the multiplies of the identity, we obtain
that the invariant Ricci operator
$$
\Ricg_{\ip}=\unc(\Ric_{\ip}-J_1\Ric_{\ip}J_1-J_2\Ric_{\ip}J_2-J_3\Ric_{\ip}J_3),
$$
satisfies $\Ricg_{\mu}|_{\ngo_i}\in\RR I$ for any $\mu\in W_h$.  By
arguing as in (\ref{multide}), we obtain that any $\mu\in W_h$ is
minimal, or equivalently, $\ip$ is a minimal compatible metric for
every $(N_{\mu},\{ J_1,J_2,J_3\})$, $\mu\in W_h$.  This proves part
(i).

In this case $K_{\gamma}=\Spe(2)$, and for each $\mu\in W$ we have
that $K_{\gamma}.\mu\cap W = \Spe(1)\times\Spe(1).\mu$.  Part (ii)
follows then from the uniqueness result in Theorem \ref{equiv1g} and
part (i).  Finally, parts (iii) and (iv) follow from an elementary
analysis of $W_h$ and $W_{ah}$ as $\Spe(1)\times\Spe(1)$-modules.
\end{proof}

We will give now an explicit continuous family of non-abelian
hypercomplex structures on $\ggo_3$, the $8$-dimensional Lie algebra
obtained as the direct sum of an abelian factor and the
$7$-dimensional quaternionic Heisenberg Lie algebra.   Such a curve
was again obtained via the variational method as in Example
\ref{m26}.

\begin{example} {\rm Consider the family defined by {\small
$$
\begin{array}{lll}
\mu_{rst}(X_1,X_2)=rZ_2, && \mu_{rst}(X_2,X_3)=(1-t)Z_4, \\
\mu_{rst}(X_1,X_3)=sZ_3, && \mu_{rst}(X_2,X_4)=-(1-s)Z_3, \\
 \mu_{rst}(X_1,X_4)=tZ_4, && \mu_{rst}(X_3,X_4)=(1-r)Z_2,
\end{array}
$$ }
which is easily seen to satisfy the integrability condition for all
$J_i$'s, though is not abelian since $\mu_{rst}(X_1,X-3)=sZ_3\ne
-(1-s)Z_3=\mu_{rst}(J_1X_1,J_1X_3)$.  The Ricci operator on the
center is given by
$$
\Ric_{\mu_{rst}}|_{\ngo_2}=\unm\left[\begin{smallmatrix} 0&&&\\ &r^2+(1-r)^2&&\\ &&s^2+(1-s)^2&\\
&&&t^2+(1-t)^2\end{smallmatrix}\right],
$$
and hence the family
$$
\left\{ (N_{\mu_{rst}},\{ J_1,J_2,J_3\},\ip):\unm\leq r\leq s\leq
t,\quad r^2+s^2+t^2-r-s-t=-\unm\right\}
$$
is pairwise non-isometric.  This gives rise a surface of pairwise
non-isomorphic non-abelian hypercomplex structures on $\ggo_3$ (see
Corollary \ref{isoiso}), since $j_{\mu_{rst}}(Z)$ is invertible for
any non-zero $Z\in\ngo_2$ which orthogonal to $Z_1$. }
\end{example}

\section{Appendix}\label{app}

We briefly recall in this appendix some features of Riemannian
geometry of left invariant metrics on nilpotent Lie groups.

Consider the vector space $\sym(\ngo)$ of symmetric real valued
bilinear forms on $\ngo$, and $\pca\subset\sym(\ngo)$ the open
convex cone of the positive definite ones (inner products), which is
naturally identified with the space of all left invariant Riemannian
metrics on $N$.  Every $\ip\in\pca$ induces a natural inner product
$g_{\ip}$ on $\sym(\ngo)$ given by
$g_{\ip}(\alpha,\beta)=\tr{A_{\alpha}A_{\beta}}$ for all
$\alpha,\beta\in\sym(\ngo)$, where $\alpha(X,Y)=g(A_{\alpha}X,Y)$.
We endow $\pca$ with the Riemannian metric $g$ given by $g_{\ip}$ on
the tangent space $\tang_{\ip}\pca=\sym(\ngo)$ for any $\ip\in\pca$.
Thus $(\pca,g)$ is isometric to the symmetric space $\Gl(n)/\Or(n)$.
E. Wilson proved that $(N,\ip)$ and $(N,\ip')$ are isometric if and
only if $\ip'=\vp.\ip:=\la\vp^{-1}\cdot,\vp^{-1}\cdot\ra$ for some
$\vp\in\Aut(\ngo)$ (see the proof of \cite[Theorem 3]{Wls}).
Therefore, although the Lie bracket $\mu$ does not play any role in
the definition of a compatible metric, it is crucial in the study of
the moduli space of compatible metrics on $(N,\gamma)$ up to
isometry.

The Ricci curvature tensor $\ricci_{\ip}$ and the Ricci operator
$\Ric_{\ip}$ of $(N,\ip)$ are given by (see \cite[7.39]{Bss}),
\begin{equation}\label{ricci}
\begin{array}{rl}
\ricci_{\ip}(X,Y)=\la\Ric_{\ip}X,Y\ra=&-\unm\displaystyle{\sum_{ij}}\la\mu(X,X_i),X_j\ra\la\mu(Y,X_i),X_j\ra \\
&+\unc\displaystyle{\sum_{ij}}\la\mu(X_i,X_j),X\ra\la\mu(X_i,X_j),Y\ra,
\end{array}
\end{equation}
for all $X,Y\in\ngo$, where $\{ X_1,...,X_n\}$ is any orthonormal
basis of $(\ngo,\ip)$.  Notice that always $\scalar(N,\ip)<0$,
unless $N$ is abelian.  It is proved in \cite{Jns} that the gradient
of the scalar curvature functional $\scalar:\pca\mapsto\RR$ is given
by
\begin{equation}\label{gradsc}
\grad(\scalar)_{\ip}=-\ricci_{\ip},
\end{equation}
and hence it follows from the properties of $\pca$ described above
that
\begin{equation}\label{ricort}
\tr{\Ric_{\ip}D}=0, \qquad\forall\; {\rm symmetric}\;
D\in\Der(\ngo),
\end{equation}
where $\Der(\ngo)$ is the Lie algebra of derivations of $\ngo$ (see
for instance \cite[(2)]{critical} for a proof of this fact).

Assume now that $\ngo$ is $2$-step nilpotent, and let $\ip$ an inner
product on $\ngo$.  Consider the orthogonal decomposition
$\ngo=\ngo_1\oplus\ngo_2$, where $\ngo_2$ is the center of $\ngo$.
Thus the Lie bracket of $\ngo$ can be viewed as a skew-symmetric
bilinear map $\mu:\ngo_1\times\ngo_1\mapsto\ngo_2$.  For each
$Z\in\ngo_2$ we define $j_{\mu}(Z):\ngo_1\mapsto\ngo_1$ by
$$
\la j_{\mu}(Z)X,Y\ra=\la\mu(X,Y),Z\ra, \qquad X,Y\in\ngo_1.
$$
$(N,\ip)$ is said to be a {\it modified H-type} Lie group if for any
non-zero $Z\in\ngo_2$
$$
j_{\mu}(Z)^2=c(Z)I \qquad \mbox{for some}\; c(Z)<0,
$$
and it is called {\it H-type} when $c(Z)=-\la Z,Z\ra$ for all
$Z\in\ngo_2$.  These metrics, introduced by A. Kaplan, play a
remarkable role in the study of Riemannian geometry on nilpotent and
solvable Lie groups (see for instance \cite{BrnTrcVnh} for further
information and \cite{modified} for the `modified' case).

If $\mu'=\vp.\mu$ for some
$\vp=(\vp_1,\vp_2)\in\Gl(\ngo_1)\times\Gl(\ngo_2)$, then it is easy
to see that
$$
j_{\mu'}(Z)=\vp_1j_{\mu}(\vp_2^tZ)\vp_1^t, \qquad \forall
Z\in\ngo_2.
$$

\end{document}